\RequirePackage{fix-cm}
\documentclass[smallextended]{svjour3}  

\let\proof\relax

\usepackage[misc]{ifsym}
\usepackage{amssymb}
\usepackage{amsthm}
\usepackage{amsmath}
\usepackage{graphicx}
\usepackage{enumerate}
\usepackage{multirow}
\usepackage{longtable}
\usepackage{booktabs}
\usepackage{float}
\usepackage{arydshln}
\usepackage{algorithmic}
\usepackage{algorithm}
\usepackage[numbers,square,sort&compress]{natbib}

\begin{document}
	\title{A New Approach to the Determination of Expert Weights in Multi-attribute Group Decision Making}
	\author{{Yuetong Liu}    \textsuperscript{1}  
	\and {Chaolang Hu}   \textsuperscript{1,~\Letter}   
	\and {Shiquan Zhang}   \textsuperscript{1,~\Letter}  
	\and {Qixiao Hu}   \textsuperscript{1} } 
	\institute{\begin{itemize}
			\item[\textsuperscript{\Letter}] {Chaolang Hu} \\
			\email{huchaolang@scu.edu.cn}
			\and \item[\textsuperscript{\Letter}] {Shiquan Zhang} \\
			\email{shiquanzhang@scu.edu.cn}
			\at 
			\item[\textsuperscript{1}] School of Mathematics, Sichuan University, Chengdu, 610065, China.
		\end{itemize}}
\date{Received: date / Accepted: date}  

   \maketitle

\begin{abstract}
This paper presents a new approach based on optimization model to determine the weights of experts in the multi-attribute group decision. Firstly, by minimizing the sum of differences between individual evaluations and the overall consistent evaluations of all experts, a new optimization model is established for determining expert weights. Then, rigorous proof of the unique existence of solution is analyzed in detail, and the sequential least squares quadratic programming algorithm is adopted to solve the optimization model. Finally, the reasonableness of the new approach is verified by numerical experiments, i.e., the smaller the difference between the individual evaluations and the overall consistent evaluations, the larger the weights assigned to the corresponding individual.

\end{abstract}

\section{Introduction}
Group decision making (GDM), as one of the main forms of decision-making activities, refers to the decision-making activities carried out by a group as the main body. Modern human decision-making activities involve a wide range of information and influencing factors, in order to achieve scientific decision-making, depending on the ability of a person is not likely to be completed, it needs to focus on the advantages of the group and the wisdom of the crowd in order to make the best decision. At the same time, the democratisation of decision-making requirements, is the basic premise and guarantee of modern scientific decision-making \cite{ref1,ref2,ref3}. Multi-criteria decision making (MCDM), which is based on the solution space of the problem under study, is categorized into two types: multi-attribute decision making (MADM) and multi-objective decision making (MODM). Multi-attribute decision making (MADM), also known as finite-alternative multi-objective decision making, refers to the decision problem of selecting the optimal alternatives or ranking them under the consideration of multiple attributes, and it is an important part of modern decision science. 

The theory and methods of multi-attribute group decision making (MAGDM) are widely used in many fields such as engineering, technology, economy, management and military, and it often has some common elements: multiple alternatives, multiple evaluation indicators, allocation of indicator weights and allocation of expert weights.

The weights of the indicators involved in evaluating a programme reflect the relative importance among attributes and are central to the accuracy of decision-making results. In recent years, scholars at home and abroad have conducted extensive research on the determination of indicator weights from different perspectives. The most common ones are AHP (Analytic Hierarchy Process) \cite{ref6,ref7,ref8}, ANP (Analytic Network Process) \cite{ref9}, TOPSIS (Technique for Order of Preference by Similarity to Ideal Solution) \cite{ref10,ref11,ref12,ref13,ref14}, ELECTRE (Elimination and Choice Expressing Reality) \cite{ref15,ref16,ref17}, VIKOR (VIseKriterijumska Optimizacija I Kompromisno Resenje) \cite {ref18}, PROMETHEE (Preference Ranking Organization Method for Enrichment Evaluations) \cite{ref19,ref20,ref21,ref22}. These methods are often simple, feasible and operational. However, the determination of expert weights is either too cumbersome or insufficiently truthful and objective, and thus has not been well addressed. In the case of complete ignorance of the information about the experts, we generally assume that the weights of each expert are equal. In practice, however, factors such as popularity, title, education, relevance to the field of study, and the number of years engaged in the relevant work, all affect the degree of the experts' contribution to the final judgement. Therefore, different experts should have different weights, coupled with the fact that different expert weights in group decision-making will produce qualitative changes in the final decision-making results, so how to rationalise expert weights has always been a hotspot for scholars to study.

Overall, there are three broad methods for setting expert weights: subjective weighting, objective weighting, and combined subjective and objective weighting.
The subjective weighting method involves determining the weights of each expert in advance, or determining the weights during the interaction of the individual experts \cite{ref23,ref24,ref25}. The weights determined by this method are not changed once they are set, i.e., they are not affected by the results of the experts' evaluations. Subjective weighting requires that the experts are well acquainted with each other and have a high level of expertise, but even so, the weights determined in this way are still highly subjective and uncertain.
Objective weighting is a method of determining the importance of information through a judgement matrix formed from the evaluation information provided by the experts, and the weights obtained change with the evaluation results \cite{ref26}. In \cite{ref27}, it addresses the supplier selection decision problem where the attribute evaluation value is a linguistic variable and the expert weights are unknown, and integrates the expert weights in terms of the hesitancy and similarity of the evaluation information in order to obtain a comprehensive evaluation matrix after aggregation, so as to realize the ranking of the suppliers through TOPSIS. \cite{ref28} introduces the degree of fault detection into the hardware FMEA method, and applies the fuzzy theory to propose an expert assignment method based on the D-S evidence theory to reasonably allocate the weights, which improves the problem of the traditional failure mode and effect analysis (FMEA), which is more subjective and the allocation of the weights is not reasonable enough.

In addition, some articles determine the expert weights from both subjective and objective perspectives in a comprehensive manner. \cite{ref29} determines the weights of attributes based on the subjective and objective value matrices of the attributes, and integrates the two weights by using the minimum entropy criterion, which fully takes into account the different judgements of all decision makers on the multi-attribute decision problem. \cite{ref30} establishes a comprehensive evaluation system, based on the fuzzy analytic hierarchy process (FAHP) method, using a combination of expert weight coefficients and entropy weighting method, and integrating subjective and objective weights according to the characteristics of the indicators, so as to ultimately determine the weights and obtain the comprehensive evaluation results.

The structure of this paper is as follows. In the first two parts, we sort out the necessity of determining expert weights in MAGDM with the current state of research, and describe the problem to be solved. In the third part, we construct an optimization problem by calculating the distance from the expert's score point to the consistent score point of the group's decision to select the expert weights which minimize the overall differences between the evaluations of each expert and the consistent results, and briefly introduce the algorithm SLSQP for solving this kind of algorithm with both equality and inequality constraints. In the fourth section, we validate the correctness and reasonableness of the proposed model by numerical experiments with the data generated by simulation. At last, we summarize the main work of this paper, and propose two directions of the future improvement of the model.

\section{Multi-attribute group decision making}
	For MAGDM, its commonalities can be written as follows:
	\begin{enumerate}[1]
	\item Multiple alternatives. Before making a group decision, the decision maker must first weigh the number of alternatives as a choice for evaluation.
	\item Multiple evaluation indicators. Before making a decision the decision maker needs to come up with relevant factors that may affect the alternative, which may be independent or interrelated. At the same time, the number of feasible indicators (attributes) needs to be identified and measured.
	\item Allocation of indicator weights. Since different indicators have different degrees of influence on the decision-making of the alternatives, the corresponding weights of different indicators also vary. In general the assignment of indicator weights is usually formalised \cite{ref4,ref5}.
	\item Allocation of expert weights. The body of GDM is a group rather than a particular individual, so the allocation of expert weights should be considered for MADM. The evaluation of diverse decision makers can have distinct effects on the final decision. Therefore, to measure these effects, weights are assigned to each decision maker based on their corresponding impact.
	\end{enumerate}

	How to objectively determine the weight of experts is a focus of current academic research. The current objective weighting method is summarized in three main ideas: the use of experts to determine the proportion of the judgement matrix scale to determine the weight of the expert, through these given values to detect the importance of the judgement information provided by the expert, i.e., the truthfulness of the information and credibility, and then determine the weights of experts \cite{ref26,ref32}; the use of the consistency ratio of judgement matrices to obtain the weights of experts which are obtained by using the eigenvectors of matrices to detect the magnitude of the contribution of the information provided by the experts \cite{ref31,ref33}; the level of judgement information provided by the experts is detected by using the value of the consistency ratio of the judgement matrix to determine the weights for each expert \cite{ref34,ref35,ref36}.

    Based on this situation, this paper establishes an optimization model for determining the weights of the experts for the evaluation results of all the experts involved in the decision-making, with the objective of minimizing the sum of the differences between the individual evaluations and the overall consistent evaluation results. 
    
    First, for multi-attribute group decision-making involving multiple alternatives, the following sets are often involved:
\begin{enumerate}[1]
	\item The set of evaluation indicators (attributes): $U=\{u_1,u_2, ... ,u_n\}$, i.e. there are $n$ attributes associated with alternatives being evaluated;
	\item The set of alternatives: $D=\{d_1,d_2,... ,d_s\}$, i.e. there are $s$ alternatives involved in the evaluation;
	\item The set of experts: $C=\{c_1,c_2,... ,c_m\}$, i.e. there are $m$ experts involved in the evaluation.
\end{enumerate}

In the actual decision-making process, it is often necessary for a number of experts to score alternatives from different perspectives (i.e., on the basis of different indicators) and to make decisions on the basis of the results. For the experts involved in the decision-making, due to the differences in their qualifications, experience and knowledge of the relevant fields, not every expert has the same authority when making decisions on the alternatives, which also indicates that it is meaningful and necessary to propose an effective method for determining the weights. It can be argued that the evaluations of the alternatives by a more experienced, knowledgeable and senior expert should be more informative, i.e. the evaluation results of that expert are more authoritative, and their importance can be demonstrated by giving them a relatively greater weight. 

In this paper, this importance can be expressed by comparing the difference between the individual evaluation and the comprehensive evaluation of the group, and the smaller the difference between the two, the smaller the conflict between the evaluation of the expert and the results of the group's consistent evaluation, and the higher the degree of consistency. Based on this idea, we propose a method for determining expert weights through an optimization model.

\section{Optimization model for expert weight determination}
\subsection{Optimization model}
Assuming that the expert $c_j$ scores the $k$-th attribute $u_k$ of alternative $d_i$ with a value of $a_{ij}^k$, the vector of scores of the expert $c_j$ for alternative $d_i$ can be written as $\boldsymbol {a_{ij}}=\begin{pmatrix}a_{ij}^1 & a_{ij}^2 & \cdots & a_{ij}^n \end{pmatrix}^T$, where $i=1,2,\dots,s$, $j=1,2,\dots,m$, and $k=1,2,\dots,n$. Then for the alternative $d_i$, its overall score is represented by the matrix $A^i$, i.e. \begin{equation}
	A^i=\begin{pmatrix}
		\boldsymbol {a_{i1}} & \boldsymbol {a_{i2}} & \cdots & \boldsymbol {a_{im}}
	\end{pmatrix}=\begin{pmatrix}
		a_{i1}^1 & \cdots & a_{im}^1\\
		\cdots & \cdots & \cdots\\
		a_{i1}^n & \cdots & a_{im}^n\end{pmatrix}. \end{equation}
For the matrix $A^i$, the $k$-th row represents the evaluation scores of the experts on the attribute $u_k$ of the alternative $d_i$, and the $j$-th column represents the evaluation scores of the experts $c_j$ on each indicator of the alternative. At the same time, assume that the expert weight vector is $\boldsymbol{w}=\begin{pmatrix}
	w_1 & w_2 & \cdots & w_m\end{pmatrix}^T$, and it is now to be determined.

The basic idea of the article is to use the results of the experts' evaluation of the multiple attributes of the alternatives to build an optimization model to objectively determine the expert weights. In addition, it should be mentioned that the experts here are generalized, all the members of the group participating in the decision-making are ``experts'' as pointed out in the article, which are essentially all the decision-makers. In actual group decision-making, not only people with professional knowledge, but often multiple groups participate in decision-making and listen to the advice of multiple groups.

$d_i$ has an ideal scoring result, which is obtained by linearly weighting the evaluation results of $m$ experts, and it reflects the consistency of the evaluation of the alternative by the expert group, which is called the ``consistent scoring point'' of $d_i$, and denoted as $\boldsymbol{b^i}$. alternatives will be evaluated from the perspective of $n$ indicators, that is to say, the ``consistent score point'' $\boldsymbol{b^i}$ of $d_i$ is an $n$-dimensional vector, that is, \begin{equation}
	\boldsymbol{b^i}=\sum_{j=1}^m w_j \boldsymbol {a_{ij}}=A^i \cdot \boldsymbol{w}.\end{equation}
A large $d_{ij}$ indicates that there is a large difference between the evaluation of expert $c_j$ and the comprehensive evaluation of the group of $d_i$, and the authority of the expert in the GDM is low. For this reason, in order to ensure the consistency of the group decision-making, this expert should be assigned a smaller weight. On this basis, the sum of the distance between each expert's evaluation of the $d_i$ and the ``consistent score point'' is recorded as $D_i$,
\begin{equation}D_i=\sum_{j=1}^m d_{ij} =\sum_{j=1}^m \Vert \boldsymbol{a_{ij}}-\boldsymbol{b^i} \Vert_2.\end{equation}.

In order to make full use of the evaluation information of all the experts on each alternative, we can integrate the data corresponding to all the alternatives to form a block matrix and determine the overall ``consistent score point'', so as to build an integrated optimization model to determine a unique set of expert weights. If we use the matrix $S$ to represent this integrated block matrix, then there are \begin{equation}S=\begin{pmatrix}
		A^1 \\ A^2 \\ \cdots \\ A^s
	\end{pmatrix}=\begin{pmatrix}
		\boldsymbol{a_{11}} & \boldsymbol{a_{12}} & \cdots & \boldsymbol{a_{1m}}\\
		\boldsymbol{a_{21}} & \boldsymbol{a_{22}} & \cdots & \boldsymbol{a_{2m}}\\
		\vdots & \vdots & \vdots & \vdots\\
		\boldsymbol{a_{s1}} & \boldsymbol{a_{s2}} & \cdots & \boldsymbol{a_{sm}}
	\end{pmatrix} _{ns \times m} \triangleq \begin{pmatrix}
		\boldsymbol {p_1} & \boldsymbol {p_2} & \cdots & \boldsymbol {p_m}
	\end{pmatrix}.\end{equation}
Define the $j$-th column vector of matrix $S$ as $\boldsymbol{p_j}$, which represents the vector of scores of experts $c_j$ on each indicator of each alternative. Among them, the ``overall consistent score point'' $\boldsymbol{b}$ is the $ns$-dimensional column vector integrating the ``consistent score point'' $\boldsymbol{b^i}$ of each alternative, which is denoted as \begin{equation}\boldsymbol{b}=\begin{pmatrix}
		\boldsymbol{b^1} \\ \boldsymbol{b^2} \\ \cdots \\ \boldsymbol{b^s}   \end{pmatrix}_{ns \times 1}=\sum_{j=1}^m w_j\boldsymbol{p_j}.\end{equation}

The distance between experts evaluation and the ``overall consistent score point'' is denoted as $s_j=\Vert \boldsymbol{p_j}-\boldsymbol{b} \Vert_2$. Based on the above idea, we build the optimization model as follows
\begin{equation}
	\begin{aligned} \label{1}
		\min_{\boldsymbol{w}} \quad & Q(\boldsymbol{w})=\sum_{j=1}^m s_j= \sum_{j=1}^m \Vert \boldsymbol{p_j}-\boldsymbol{b} \Vert_2\\
		\mbox{s.t.}\quad
		&\sum_{j=1}^m w_j =1,\\
		&0 \leq w_j \leq 1 ,j=1,2,\dots,m.
	\end{aligned}
\end{equation}
In the model, $Q(\boldsymbol{w})$ is the objective function, which represents the sum of the distances from $m$ expert evaluation points to the ``overall consistent score point'' of the decision, and we want to minimize this objective function, i.e., to minimize the sum of the differences between the experts' evaluations and the overall evaluations, in order to determine the final expert weights. The expert weights obtained on the basis of the objective function make the expert evaluations more consistent and the differences between expert evaluations as small as possible. The optimal solution to the model is the finalized ideal expert weights, which can be used by substituting the results into various evaluation models.

On one hand, the evaluation of each expert should be consistent, and if the difference between individual evaluations is too large, it will easily lead to bias in the final comprehensive evaluation results.For this reason, we hope that the distance between the evaluation points of each expert and the ``overall consistent score point'' is as small as possible, so that expert weights obtained will make the consistency of GDM relatively high.At the same time, for individual experts, when the distance between its evaluation and the overall consistent evaluation is close, we hope that the expert has a high weight, that is to say, if the expert's score is more reliable and authoritative, we will give it a greater weight. On the other hand, we default that the experts participating in the scoring are all specialized, so in order to ensure the rationality and professionalism of the overall evaluation, the gap between the decision integrating the evaluation of all experts and  each expert should not be too large. This also indicates that the expert weights obtained when the objective function we set is minimized make the consistency between the integrated decision and the evaluation of each expert the highest, then the professionalism and credibility of the final evaluation results can be guaranteed.

In addition, the objective method of determining expert weights in the existing literature and practical operation is often based on the judgment matrix given by experts, which requires experts to measure the importance of each evaluation indicator, and this operation is subjective. Coupled with the relative complexity of this process, it is seldom used in real problems, while allowing experts to evaluate directly based on the indicators is more feasible and easy to operate, compared to the existing methods are more in line with the practical application of the scene. Therefore, the method of establishing the optimization model \eqref{1} to determine the expert weights is reasonable and practical value.

\subsection{Existence and uniqueness of model solutions}
For the above optimization model, we first introduce some theories related to convex optimization that provide assistance in determining the existence and uniqueness of the solution to the optimization model \eqref{1}.
\begin{definition}[Convex function \cite{ref37}] \label{def1}
	Let $S \subset \mathbb{R}^n$ be a non-empty convex set and $f$ be a function defined on $S$. Call the function $f$ a convex function on $S$ if for any $\boldsymbol{x_1} \in S$, $\boldsymbol{x_2} \in S$, $\lambda \in (0,1)$, there are
	\begin{equation}
		f(\lambda \boldsymbol{x_1}+(1-\lambda)\boldsymbol{x_2}) \leq \lambda f(\boldsymbol{x_1})+(1-\lambda)f(\boldsymbol{x_2}).
	\end{equation}
	Call $f$ a strictly convex function on $S$ if the strict inequality sign in the above equation holds when $\boldsymbol{x_1} \neq\boldsymbol{x_2}$.
\end{definition}

\begin{definition}[Convex optimization \cite{ref37}] \label{def2}
	Call the constrained problem
	\begin{equation}
		\begin{aligned} \label{2}
			\min \quad & f(\boldsymbol{x}) , \boldsymbol{x} \in \mathbb{R}^n   & \\
			\mbox{s.t.}\quad
			&c_s(\boldsymbol{x}) \leq 0 ,s=1,2,\dots,k, \\
			&d_t(\boldsymbol{x}) = 0 ,t=1,2,\dots,l,
		\end{aligned}
	\end{equation}
	a convex optimization problem if both the objective function $f(\boldsymbol{x})$ and the constraint function $c_s (\boldsymbol{x})(s=1,2,\dots,k)$ are convex and $d_t (\boldsymbol{x})(t=1,2,\dots,l)$ is linear.
\end{definition}

\begin{lemma}[Existence of convex optimization solutions \cite{ref38}] \label{lem1}
	Consider the convex optimization problem \eqref{2} and let $D$ be the feasible region of the problem, i.e.
	\begin{equation}
		D=\{ \boldsymbol{x} \vert c_s (\boldsymbol{x}) \leq 0,s=1,2,\dots,k;d_t (\boldsymbol{x}) = 0,t=1,2,\dots,l;\boldsymbol{x} \in \mathbb{R}^n \}.
	\end{equation}
	If the feasible region $D$ is a nonempty compact set and $f$ is continuous on $D$, the optimal solution $\boldsymbol{x^*}$ must exist.
\end{lemma}

\begin{lemma}[Uniqueness of convex optimization solutions \cite{ref39}] \label{lem2}
	Consider the convex optimization problem \eqref{2} and let $D$ be the feasible region of the problem, then
	\begin{enumerate}[i)]
		\item If the problem has a local optimal solution $\boldsymbol{x^*}$, then $\boldsymbol{x^*}$ is the global optimal solution to the problem;
		\item The set consisting of the global optimal solutions to the problem is a convex set;
		\item If the problem has a local optimal solution $\boldsymbol{x^*}$ and $f(\boldsymbol{x})$ is a strictly convex function on $D$, then $\boldsymbol{x^*}$ is the unique global optimal solution to the problem.
	\end{enumerate}
\end{lemma}

\begin{theorem}\label{the1}
	The objective function $Q(\boldsymbol{w})$ is strictly convex when the columns of the matrix $S$ are full rank.
\end{theorem}
\proof
	$Q(\boldsymbol{w})=\sum_{j=1}^m \Vert \boldsymbol{p_j}-\boldsymbol{b} \Vert_2=\sum_{j=1}^m \Vert \boldsymbol{p_j}-\sum_{i=1}^m w_i\boldsymbol{p_i} \Vert_2$, suppose $\boldsymbol{w}$ and $\boldsymbol{w'}$ are two sets of weights that satisfy the constraints and $\boldsymbol{w}\neq\boldsymbol{w'}$, $\lambda \in (0,1)$. Then we can get
	\begin{equation} \label{4}
		\begin{aligned}
			Q(\lambda \boldsymbol{w}+(1-\lambda)\boldsymbol{w'}) =\sum_{j=1}^m \Vert\boldsymbol{p_j}-\sum_{i=1}^m [\lambda w_i+(1-\lambda)w'_i]\boldsymbol{p_i} \Vert_2
		\end{aligned}
	\end{equation}
	and
	\begin{equation} \label{5}
		\begin{aligned}
			\lambda Q(\boldsymbol{w})+(1-\lambda)Q(\boldsymbol{w'})
			=\lambda \sum_{j=1}^m\Vert \boldsymbol{p_j}-\sum_{i=1}^m w_i \boldsymbol {p_i} \Vert_2 + (1-\lambda)\sum_{j=1}^m\Vert \boldsymbol{p_j}-\sum_{i=1}^m w'_i \boldsymbol {p_i} \Vert_2.
		\end{aligned}
	\end{equation}
	Let $F_j=\lambda \Vert \boldsymbol{p_j}-\sum_{i=1}^m w_i \boldsymbol {p_i} \Vert_2 + (1-\lambda)\Vert \boldsymbol{p_j}-\sum_{i=1}^m w'_i \boldsymbol {p_i} \Vert_2$ as well as $H_j=\Vert\boldsymbol{p_j}-\sum_{i=1}^m [\lambda w_i+(1-\lambda)w'_i]\boldsymbol{p_i} \Vert_2$, remember that $T(\lambda;\boldsymbol{w},\boldsymbol{w'})=\lambda Q(\boldsymbol{w})+(1-\lambda)Q(\boldsymbol{w'})-Q(\lambda \boldsymbol{w}+(1-\lambda)\boldsymbol{w'})$, so we have
	\begin{equation} \label{6}
		\begin{aligned}
			&T(\lambda;\boldsymbol{w},\boldsymbol{w'})\\
			=&\lambda Q(\boldsymbol{w})+(1-\lambda)Q(\boldsymbol{w'})-Q(\lambda \boldsymbol{w}+(1-\lambda)\boldsymbol{w'}) \\
			=&\lambda \sum_{j=1}^m\Vert \boldsymbol{p_j}-\sum_{i=1}^m w_i \boldsymbol {p_i} \Vert_2 + (1-\lambda)\sum_{j=1}^m\Vert \boldsymbol{p_j}-\sum_{i=1}^m w'_i \boldsymbol {p_i} \Vert_2-\sum_{j=1}^m \Vert\boldsymbol{p_j}-\sum_{i=1}^m [\lambda w_i+(1-\lambda)w'_i]\boldsymbol{p_i} \Vert_2\\
			=&\sum_{j=1}^m\{\lambda \Vert \boldsymbol{p_j}-\sum_{i=1}^m w_i \boldsymbol {p_i} \Vert_2 + (1-\lambda)\Vert \boldsymbol{p_j}-\sum_{i=1}^m w'_i \boldsymbol {p_i} \Vert_2- \Vert\boldsymbol{p_j}-\sum_{i=1}^m [\lambda w_i+(1-\lambda)w'_i]\boldsymbol{p_i} \Vert_2\}\\
			\triangleq & \sum_{j=1}^m (F_j-H_j).
		\end{aligned}
	\end{equation}
	
	Since $\lambda \in (0,1)$ is positive, and the L2 norm of any vector must be non-negative, both $F_j$ and $H_j$ are non-negative. Square both, and take the difference as
	\begin{equation} \label{7}
		\begin{aligned}
			&F_j^2-H_j^2\\
			=&[\lambda \Vert \boldsymbol{p_j}-\sum_{i=1}^m w_i \boldsymbol {p_i} \Vert_2 + (1-\lambda)\Vert \boldsymbol{p_j}-\sum_{i=1}^m w'_i \boldsymbol {p_i} \Vert_2]^2- \Vert\boldsymbol{p_j}-\sum_{i=1}^m [\lambda w_i+(1-\lambda)w'_i]\boldsymbol{p_i} \Vert_2^2 \\
			=&\lambda^2\Vert \boldsymbol{p_j}-\sum_{i=1}^m w_i \boldsymbol {p_i} \Vert_2^2+(1-\lambda)^2\Vert \boldsymbol{p_j}-\sum_{i=1}^m w'_i \boldsymbol {p_i} \Vert_2^2+2\lambda(1-\lambda) \Vert \boldsymbol{p_j}-\sum_{i=1}^m w_i \boldsymbol {p_i} \Vert_2\Vert \boldsymbol{p_j}-\sum_{i=1}^m w'_i \boldsymbol {p_i} \Vert_2 \\
			&-\Vert\lambda(\boldsymbol{p_j}-\sum_{i=1}^mw_i\boldsymbol{p_i})+(1-\lambda)(\boldsymbol{p_j}-\sum_{i=1}^m w'_i\boldsymbol{p_i} )\Vert_2^2 \\
			=&\lambda^2\Vert \boldsymbol{p_j}-\sum_{i=1}^m w_i \boldsymbol {p_i} \Vert_2^2+(1-\lambda)^2\Vert \boldsymbol{p_j}-\sum_{i=1}^m w'_i \boldsymbol {p_i} \Vert_2^2+2\lambda(1-\lambda) \Vert \boldsymbol{p_j}-\sum_{i=1}^m w_i \boldsymbol {p_i} \Vert_2\Vert \boldsymbol{p_j}-\sum_{i=1}^m w'_i \boldsymbol {p_i} \Vert_2 \\
			&-\lambda^2\Vert\boldsymbol{p_j}-\sum_{i=1}^mw_i\boldsymbol{p_i}\Vert_2^2-(1-\lambda)^2\Vert\boldsymbol{p_j}-\sum_{i=1}^m w'_i\boldsymbol{p_i}\Vert_2^2-2\lambda(1-\lambda)(\boldsymbol{p_j}-\sum_{i=1}^mw_i\boldsymbol{p_i})^T(\boldsymbol{p_j}-\sum_{i=1}^m w'_i\boldsymbol{p_i}) \\
			=&2\lambda(1-\lambda)[ \Vert \boldsymbol{p_j}-\sum_{i=1}^m w_i \boldsymbol {p_i} \Vert_2\Vert \boldsymbol{p_j}-\sum_{i=1}^m w'_i \boldsymbol {p_i} \Vert_2-(\boldsymbol{p_j}-\sum_{i=1}^mw_i\boldsymbol{p_i})^T(\boldsymbol{p_j}-\sum_{i=1}^m w'_i\boldsymbol{p_i})].
		\end{aligned}
	\end{equation}
	Since $\lambda \in (0,1)$, then $\lambda(1-\lambda)$ must be a real number greater than 0. For the Cauchy-Schwarz inequality
	\begin{equation} \label{8}
		\begin{aligned}
			\boldsymbol{x}^T\boldsymbol{y}=|(\boldsymbol{x},\boldsymbol{y})|\leq\Vert \boldsymbol{x}\Vert_2\Vert \boldsymbol{y} \Vert_2,
		\end{aligned}
	\end{equation}
	the equality sign holds if and only if $\boldsymbol{x}$ and $\boldsymbol{y}$ are linearly related. From this inequality it follows that
	\begin{equation} \label{9}
		\begin{aligned}
			\Vert \boldsymbol{p_j}-\sum_{i=1}^m w_i \boldsymbol {p_i} \Vert_2\Vert \boldsymbol{p_j}-\sum_{i=1}^m w'_i \boldsymbol {p_i} \Vert_2\geq(\boldsymbol{p_j}-\sum_{i=1}^mw_i\boldsymbol{p_i})^T(\boldsymbol{p_j}-\sum_{i=1}^m w'_i\boldsymbol{p_i}),
		\end{aligned}
	\end{equation}
	from this it follows that $F_j^2-H_j^2\geq0$ always holds. As $F_j$ and $H_j$ are both non-negative, we get $F_j\geq H_j$, which is $T(\lambda;\boldsymbol{w},\boldsymbol{w'})\geq0$, and then it follows that the objective function $Q(\boldsymbol{w})$ must be a convex function.
	
	Since $\boldsymbol{w}\neq\boldsymbol{w'}$, for the inequality sign of \eqref{6} to hold strictly, that is, $F_j- H_j=0$ does not hold simultaneously for $j=1,2,\dots,m$. To prove this, we use proof by contradiction, i.e., assume that they hold at the same time. Then for any $j$ it is satisfied that the vector $\boldsymbol{p_j}-\sum_{i=1}^mw_i\boldsymbol{p_i}$ is linearly related to the vector $\boldsymbol{p_j}-\sum_{i=1}^m w'_i\boldsymbol{p_i}$, which means that
	\begin{equation} \label{10}
		\begin{aligned}
			\boldsymbol{p_j}-\sum_{i=1}^mw_i\boldsymbol{p_i}=k_j(\boldsymbol{p_j}-\sum_{i=1}^m w'_i\boldsymbol{p_i}),\ \forall j.
		\end{aligned}
	\end{equation} Since the matrix columns are assumed to be full rank, which means that the matrix column vectors are linearly independent, then the coefficients of the vector $p_j$ are all zero, giving us
	\begin{equation} \label{11} \left\{
		\begin{aligned}
			&-w_i=-k_jw'_i,&i\neq j,\\
			&1-w_i=k_j(1-w'_i),&i=j.
		\end{aligned}\right.
	\end{equation}
	For $j=1,2,\dots,m,i\neq j$, we have
	\begin{equation} \label{12} \left\{
		\begin{aligned}
			&w_i=k_1w'_i,&i\neq 1,\\
			&w_i=k_2w'_i,&i\neq 2,\\
			&\cdots\\
			&w_i=k_mw'_i,&i\neq m.
		\end{aligned}\right.
	\end{equation}
	Also, for $i=j=1,2,...,m$, we can get
	\begin{equation} \label{13} \left\{
		\begin{aligned}
			&1-w_1=k_1(1-w'_1),\\
			&1-w_2=k_2(1-w'_2),\\
			&\cdots\\
			&1-w_m=k_m(1-w'_m).
		\end{aligned}\right.
	\end{equation}
	From \eqref{12}, we can deduce that $k_1=k_2=... =k_m\triangleq C$ and $w_i=Cw'_i$. Combining with \eqref{13} gives $C=1$, i.e., $\boldsymbol{p_j}-\sum_{i=1}^mw_i\boldsymbol{p_i}=\boldsymbol{p_j}-\sum_{i=1}^m w'_i\boldsymbol{p_i}$, while $w_i =w'_i$, holds for any $i=1,2,\dots,m$, which contradicts the hypothesis, thus showing that $F_j-H_j$ cannot be taken equal at the same time.
	
	From the result of proof by contradiction, we can know that $F_j - H_j=0$ cannot be taken equal at the same time, that is to say, the inequality sign of \eqref{6} is strictly valid, and thus we conclude that the objective function $Q(\boldsymbol{w})$ is strictly convex when the matrix column is full of rank.\qed

\begin{theorem}\label{the2}
	If the matrix $S$ columns are full rank, then the solution to the optimization problem \eqref{1} exists and it is unique.
\end{theorem}
\proof
	To facilitate the analysis of the problem, we first represent \eqref{1} in the following standard format
	\begin{equation}
		\begin{aligned} \label{3}
			\min_{\boldsymbol{w}} \quad & Q(\boldsymbol{w})= \sum_{j=1}^m \Vert \boldsymbol{p_j}-\boldsymbol{b} \Vert_2\\
			\mbox{s.t.}\quad
			&w_j - 1 \leq 0 ,j=1,2,\dots,m,\\
			&- w_j \leq 0 ,j=1,2,\dots,m,\\
			&\sum_{j=1}^m w_j -1 = 0,
		\end{aligned}
	\end{equation}where $w_j-1 $ and $- w_j $ are both convex functions and the equality constraint $\sum_{j=1}^m w_j -1 = 0$ is linear. Since $Q(\boldsymbol{w})$ is convex, by definition \ref{def2} it follows that \eqref{1} is a convex optimization problem. For this reason by Lemma \ref{lem1} we conclude that a local optimal solution to \eqref{1} must exist.
	By the conclusion of Theorem \ref{the1} we know that in the case of matrix columns with full rank, $Q(\boldsymbol{w})$ is not only a convex function, but also strictly convex. Then by Lemma \ref{lem2} we know that the local optimal solution of this optimization model is the global optimal solution, and the solution exists and is unique, so the theorem can be proved.\qed

In addition to this, if there are some additional qualifications and constraints on the alternatives before they are chosen, then we can also implement constraints on the expert weights by writing them into the constraints of the optimization problem. In the real engineering field, many times there are other constraints on the problem, for example, before making the decision, we already know that expert $a$ is more specialized than expert $b$ in evaluating the solution, so the importance of expert $a$ evaluation should be higher than that of expert $b$. Accordingly we add a new constraint $w_a>w_b$ for \eqref{1} to ensure that the expert weights obtained by solving \eqref{1} meet the requirements. The addition of constraints is flexible, as long as the addition of constraints still ensures that \eqref{1} is a convex optimization model, then we will be able to obtain the unique solution under this model, i.e., the expert weights under this decision model.

Matrix $S$ is a matrix with $ns$ rows and $m$ columns, and its rows are equal to the product of the number of indicators $n$ and the number of alternatives $s$. In practice, we often want the number of decision makers to increase and the number of alternatives and indicators to be as large as possible. From the matrix point of view, we would also like to see the number of experts and the indicators involved in the evaluation of the decision-making process to be as large as possible, so that the number of rows of the matrix $S$ increases, and thus the likelihood that the columns of the matrix will be full-ranked increases.

\subsection{SLSQP: a algorithm for solving optimization models}
The optimization model proposed above is a nonlinear optimization problem with both equality and inequality constraints, and for this type of optimization problem, we often choose to use the SQP (Sequential Quadratic Programming) algorithm to solve it.The SQP method is one of the most efficient computational methods for solving general nonlinear programming problems. The core idea of the algorithm is to construct quadratic programming subproblems in each iteration and to approximate the optimal solution step by step by solving the subproblems \cite{ref40}.
For general nonlinear problems\begin{equation}\begin{aligned}\label{21}
		\min_{\boldsymbol{x}\in \mathbb{R}^n} \quad & f(\boldsymbol{x})\\
		\mbox{s.t.}\quad
		&g_j(\boldsymbol{x})=0,j=1,2,\dots,m_e,\\
		&g_j(\boldsymbol{x})\geq 0,j=m_e+1,m_e+2,\dots,m,
\end{aligned}\end{equation}
where $f:\mathbb{R}^n \rightarrow \mathbb{R}^1$ and $\boldsymbol{g}:\mathbb{R}^n \rightarrow \mathbb{R}^m$ are continuously differentiable functions without any special structure. Now find the local minimum solution to the above nonlinear programming problem. We can solve \eqref{21} by iteration, assuming that given the initial vector $\boldsymbol{x_0}$, the $k+1$-th iteration can be obtained by the following equation \begin{equation}
	\boldsymbol{x_{k+1}}=\boldsymbol{x_k}+\alpha_k \boldsymbol{d_k}, \end{equation}
where $\boldsymbol{d_k}$ is the search direction and $\alpha_k$ is the step size.

The search direction $\boldsymbol{d_k}$ is determined by the quadratic programming subproblem. This subproblem consists of the quadratic approximation of the Lagrange function $\mathcal{L}(\boldsymbol{x},\boldsymbol{\lambda})=f(\boldsymbol{x})-\boldsymbol{\lambda}^T\boldsymbol{g}(\boldsymbol{x})$ and the bounded linear approximation, usually written in the following form
\begin{equation}\begin{aligned}\label{23}
		\min_{\boldsymbol{d}} \quad & \nabla f_k^T\boldsymbol{d}+\frac{1}{2}\boldsymbol{d}^TB_k\boldsymbol{d}\\
		\mbox{s.t.}\quad
		&\nabla g_j(\boldsymbol{x_k})^T\boldsymbol{d}+g_j(\boldsymbol{x_k})=0,j=1,2,\dots,m_e,\\
		&\nabla g_j(\boldsymbol{x_k})^T\boldsymbol{d}+g_j(\boldsymbol{x_k})\geq 0,j=m_e+1,m_e+2,\dots,m.
\end{aligned}\end{equation}

The 1-dimensional non-differentiable exact penalty function \begin{equation}\phi_{\boldsymbol{\rho}}(\boldsymbol{x})=f(\boldsymbol{x})+\sum_{i=1}^{m_e}\rho_j|g_j(\boldsymbol{x})|+\sum_{i=m_e+1}^m\rho_j|g_j(\boldsymbol{x})|_{-}\label{24}\end{equation}
can be used as a merit function\begin{equation}\varphi(\boldsymbol{\alpha})=\phi_{\boldsymbol{\rho}}(\boldsymbol{x_k}+\alpha \boldsymbol{d_k}),\label{26}\end{equation}
where $|g_j(\boldsymbol{x})|_{-}=|min(0,g_j(\boldsymbol{x}))|$, $\boldsymbol{x_k}$ and $\boldsymbol{d_k}$ are fixed. The update of the penalty parameter is written as $\rho_j=max\{\frac{1}{2}(\rho_j^{-}+|\mu_j|),|\mu_j|\},j=1,2,\dots,m$. Minimizing the merit function thus allows to solve for the corresponding step size $\alpha_k$. In the merit function, $\boldsymbol{\mu}$ denotes the Lagrange multiplier of the quadratic subproblem \eqref{23}, and $\rho_j^{-}$ is the $j$-th penalty parameter of the previous iteration ($\boldsymbol{\rho^0}=0$). In addition to this, if a differentiable augmented Lagrangian function is used instead of \eqref{24}, the difficulties that may be encountered in line search of non-differentiable merit functions can be overcome.

In practice, the matrix $B_k$ is not computed at every iteration to improve computational efficiency when applying SQP, for which we update the matrix $B_k$ with the following correction to ensure that it remains positive definite under any initial estimation.
\begin{equation}B_{k+1}=B_k+\frac{\boldsymbol{q_k}\boldsymbol{q_k}^T}{\boldsymbol{q_k}^T\boldsymbol{s_k}}-\frac{B_k\boldsymbol{s_k}\boldsymbol{s_k}^TB_k}{\boldsymbol{s_k}^TB_k\boldsymbol{s_k}},\label{27}\end{equation}
\begin{equation}\boldsymbol{s_k}=\boldsymbol{x_{k+1}}-\boldsymbol{x_k}=\alpha_k\boldsymbol{d_k},\end{equation}
\begin{equation}\boldsymbol{q_k}=\theta_k\boldsymbol{\eta_k}+(1-\theta_k)B_k\boldsymbol{s_k},\end{equation}
where $\boldsymbol{\eta_k}$ is the difference between the gradient of the Lagrangian function,
\begin{equation}\boldsymbol{\eta_k}=\nabla\mathcal{L}(\boldsymbol{x_{k+1}},\boldsymbol{\lambda_k})-\nabla\mathcal{L}(\boldsymbol{x_k},\boldsymbol{\lambda_k}).\end{equation}
The choice of $\theta_k$ is as follows
\begin{equation}
	\theta_k=\left\{\begin{aligned}
		&1, \boldsymbol{s_k}^T\boldsymbol{\eta_k}\geq 0.2\boldsymbol{s_k}^TB_k\boldsymbol{s_k},\\
		&\frac{0.8\boldsymbol{s_k}^TB_k\boldsymbol{s_k}}{\boldsymbol{s_k}^TB_k\boldsymbol{s_k}-\boldsymbol{s_k}^T\boldsymbol{\eta_k}}, otherwise.\end{aligned}\right.
	\label{28}\end{equation}
It ensures that the condition $\boldsymbol{s_k}^T\boldsymbol{q_k}\geq 0.2\boldsymbol{s_k}^TB_k\boldsymbol{s_k}$ holds, that is, $B_{k+1}$  is positively definite within the linear manifold defined by $\boldsymbol{x _{k+1}}$ within the linear manifold defined by the tangent plane to active constraints at $\boldsymbol{x_{k+1}}$.

The optimization problem is often solved using Python by the SQP method using the SLSQP (Sequential Least Squares Quadratic Programming). SLSQP is solved by Cholesky decomposition of the matrix $B$, i.e., letting $B=LDL^T$. Thus, the problem is equivalently transformed into an optimization problem where the objective function is of least squares form like \eqref{25} \cite{ref40,ref41}.
\begin{equation}\begin{aligned}\label{25}
		\min_{d} \quad & \Vert D_k^{1/2}L_k^T\boldsymbol{d}+ D_k^{-1/2}L_k^{-1}\nabla f_k\Vert_2\\
		\mbox{s.t.}\quad
		&\nabla g_j(\boldsymbol{x_k})^T\boldsymbol{d}+g_j(\boldsymbol{x_k})=0,j=1,2,...,m_e,\\
		&\nabla g_j(\boldsymbol{x_k})^T\boldsymbol{d}+g_j(\boldsymbol{x_k})\geq 0,j=m_e+1,m_e+2,...,m.
\end{aligned}\end{equation}

The algorithm of SLSQP is summarized below:
\begin{algorithm}[]
	\caption{Algorithm of SLSQP}\label{alg:A}
	\begin{algorithmic}[1]
		\STATE  $\boldsymbol{x_0}$, $B_0=I$, $k=0$
		\WHILE{convergence is not arise}
		\STATE compute$\nabla f_k$, $\nabla g_i(\boldsymbol{x_k})$, $g_i(\boldsymbol{x_k})$
		\STATE compute $L_k$, $D_k$ by $B_k=L_kD_kL_k^T$
		\STATE solve \eqref{25} to obtain the solution $\boldsymbol{d_k}$ and the corresponding Lagrange multiplier $\boldsymbol{\mu_k}$
		\STATE solve \eqref{26} to obtain the step size $\alpha_k$
		\STATE set $\boldsymbol{x_{k+1}}=\boldsymbol{x_k}+\alpha_k\boldsymbol{d_k}$
		\STATE compute $B_{k+1}$ by \eqref{27}-\eqref{28}
		\STATE increase $k$ by 1
		\ENDWHILE
		\RETURN $\boldsymbol{x_k}$ and $\boldsymbol{\mu_k}$ as $\boldsymbol{x^*}$ and $\boldsymbol{\mu^*}$ respectively
\end{algorithmic}\end{algorithm}

\section{Numerical Experiment}
Based on the above proposed model, we randomly generated several sets of data and conducted numerical experiments on the computer, and obtained results with certain regularity. In the following, we introduce and analyze the proposed method by taking the problem consisting of this set of data in Table \ref{tab1} as an example.

\begin{table}[]
	\centering
	\caption{Expert evaluations under each indicator for alternatives.}
		\resizebox{\textwidth}{!}{
	\begin{tabular}{|c|c|c|c|c|c|c|c|c|}
		\hline
		\textbf{Alternatives} 	& \textbf{Indicators}     & \textbf{$\boldsymbol{c_1}$} & \textbf{$\boldsymbol{c_2}$} & \textbf{$\boldsymbol{c_3}$} & \textbf{$\boldsymbol{c_4}$} & \textbf{$\boldsymbol{c_5}$} & \textbf{$\boldsymbol{c_6}$} & \textbf{$\boldsymbol{c_7}$} \\ \hline
		\multirow{6}{*}{\textbf{$\boldsymbol{d_1}$}} & \textbf{$\boldsymbol{u_1}$} & 55           & 86           & 61           & 75           & 44           & 47           & 95           \\ \cline{2-9} 
		& \textbf{$\boldsymbol{u_2}$} & 64           & 76           & 48           & 61           & 72           & 40           & 87           \\ \cline{2-9} 
		& \textbf{$\boldsymbol{u_3}$} & 48           & 44           & 76           & 80           & 81           & 68           & 93           \\ \cline{2-9} 
		& \textbf{$\boldsymbol{u_4}$} & 78           & 74           & 80           & 51           & 41           & 42           & 48           \\ \cline{2-9} 
		& \textbf{$\boldsymbol{u_5}$} & 52           & 45           & 67           & 62           & 83           & 53           & 50           \\ \cline{2-9} 
		& \textbf{$\boldsymbol{u_6}$} & 72           & 72           & 86           & 53           & 50           & 70           & 91           \\ \hline
		\multirow{6}{*}{\textbf{$\boldsymbol{d_2}$}} & \textbf{$\boldsymbol{u_1}$} & 97           & 71           & 63           & 70           & 81           & 54           & 62           \\ \cline{2-9} 
		& \textbf{$\boldsymbol{u_2}$} & 94           & 77           & 86           & 94           & 64           & 85           & 95           \\ \cline{2-9} 
		& \textbf{$\boldsymbol{u_3}$} & 65           & 92           & 63           & 54           & 92           & 65           & 86           \\ \cline{2-9} 
		& \textbf{$\boldsymbol{u_4}$} & 77           & 97           & 52           & 99           & 86           & 66           & 69           \\ \cline{2-9} 
		& \textbf{$\boldsymbol{u_5}$} & 99           & 40           & 77           & 76           & 86           & 92           & 78           \\ \cline{2-9} 
		& \textbf{$\boldsymbol{u_6}$} & 52           & 93           & 87           & 90           & 50           & 98           & 59           \\ \hline
		\multirow{6}{*}{\textbf{$\boldsymbol{d_3}$}} & \textbf{$\boldsymbol{u_1}$} & 89           & 59           & 74           & 84           & 96           & 94           & 70           \\ \cline{2-9} 
		& \textbf{$\boldsymbol{u_2}$} & 90           & 77           & 46           & 80           & 93           & 95           & 72           \\ \cline{2-9} 
		& \textbf{$\boldsymbol{u_3}$} & 45           & 53           & 68           & 77           & 54           & 90           & 63           \\ \cline{2-9} 
		& \textbf{$\boldsymbol{u_4}$} & 43           & 88           & 99           & 86           & 68           & 71           & 95           \\ \cline{2-9} 
		& \textbf{$\boldsymbol{u_5}$} & 44           & 63           & 68           & 88           & 67           & 63           & 45           \\ \cline{2-9} 
		& \textbf{$\boldsymbol{u_6}$} & 84           & 68           & 54           & 51           & 47           & 52           & 55           \\ \hline
		\multirow{6}{*}{\textbf{$\boldsymbol{d_4}$}} & \textbf{$\boldsymbol{u_1}$} & 87           & 42           & 42           & 52           & 72           & 98           & 72           \\ \cline{2-9} 
		& \textbf{$\boldsymbol{u_2}$} & 84           & 71           & 64           & 95           & 90           & 44           & 87           \\ \cline{2-9} 
		& \textbf{$\boldsymbol{u_3}$} & 83           & 73           & 59           & 51           & 75           & 97           & 66           \\ \cline{2-9} 
		& \textbf{$\boldsymbol{u_4}$} & 45           & 91           & 87           & 76           & 93           & 97           & 48           \\ \cline{2-9} 
		& \textbf{$\boldsymbol{u_5}$} & 74           & 83           & 87           & 69           & 96           & 67           & 61           \\ \cline{2-9} 
		& \textbf{$\boldsymbol{u_6}$} & 57           & 51           & 51           & 80           & 64           & 88           & 50           \\ \hline
		\multirow{6}{*}{\textbf{$\boldsymbol{d_5}$}} & \textbf{$\boldsymbol{u_1}$} & 42           & 79           & 93           & 70           & 51           & 71           & 82           \\ \cline{2-9} 
		& \textbf{$\boldsymbol{u_2}$} & 46           & 93           & 78           & 49           & 86           & 86           & 55           \\ \cline{2-9} 
		& \textbf{$\boldsymbol{u_3}$} & 93           & 51           & 71           & 96           & 44           & 75           & 46           \\ \cline{2-9} 
		& \textbf{$\boldsymbol{u_4}$} & 80           & 73           & 77           & 70           & 45           & 58           & 70           \\ \cline{2-9} 
		& \textbf{$\boldsymbol{u_5}$} & 54           & 47           & 80           & 68           & 43           & 75           & 86           \\ \cline{2-9} 
		& \textbf{$\boldsymbol{u_6}$} & 78           & 63           & 65           & 57           & 48           & 71           & 46           \\ \hline
	\end{tabular}}\label{tab1}\end{table}

We select 7 experts $C=\{c_1,c_2,... ,c_7\}$, from 6 evaluation indicators $U=\{u_1,u_2,... ,u_6\}$ on 5 alternatives $D=\{d_1,d_2,... ,d_5\}$ for evaluation, and make decisions based on these scoring data. According to the way the evaluation matrix of each alternative is built to get $A^i$, put them together so that the matrix $S$ can be built, then the optimization model is built, and assuming that expert weights are $\boldsymbol{w}=\begin{pmatrix}
	w_1, & w_2, & w_3, & w_4, & w_5, & w_6, & w_7 \end{pmatrix}^T$. In solving the optimization problem, the error is chosen to be $\epsilon=10^{-12}$, and the initial weight vector is the mean weight case, $\boldsymbol{w_0}=\begin{pmatrix}1/7, &1/7, & 1/7, & 1/7, & 1/7, & 1/7, & 1/7 \end{pmatrix}^T$.

\begin{figure}[H]
	\centering
	\includegraphics[width=0.8\textwidth]{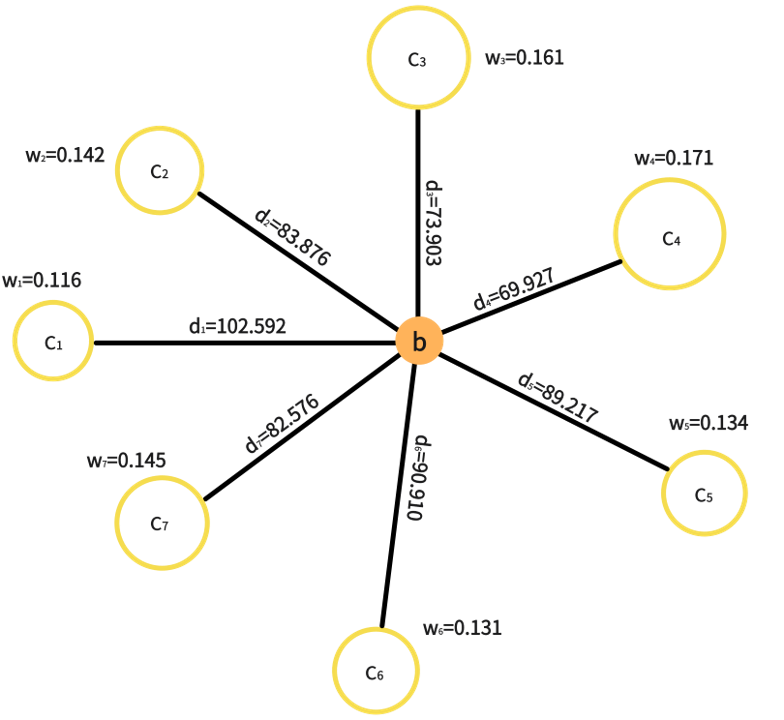}
	\caption{Weight of each expert and the distance to ``overall consistent score points''.}
	\label{Fig.1}
\end{figure}

The above problem is solved using SLSQP and the results obtained are shown in the following Table \ref{tab2} and Figure \ref{Fig.1}.
\begin{table}[H]
	\centering
	\caption{Expert weights and corresponding distances to ``overall consistent score points''.}
		\resizebox{\textwidth}{!}{
	\begin{tabular}{|c|c|c|c|c|c|c|c|}
		\hline
		& \textbf{$\boldsymbol{c_1}$} & \textbf{$\boldsymbol{c_2}$} & \textbf{$\boldsymbol{c_3}$} & \textbf{$\boldsymbol{c_4}$} & \textbf{$\boldsymbol{c_5}$} & \textbf{$\boldsymbol{c_6}$} & \textbf{$\boldsymbol{c_7}$}\\ \hline
		\textbf{Expert weight}      &0.116 &0.142 &0.161 &0.171 &0.134 &0.131
		&0.145  \\ \hline
		Descending order of weights     &7     &4       & 2       & 1     &5     & 6   & 3   \\ \hline
		\textbf{Distances} &102.592 &83.876 &73.903 &69.927 &89.217 &90.910 &82.576  \\ \hline
		Ascending order of distances     &7     &4       & 2       & 1     &5     & 6   & 3  \\ \hline
	\end{tabular}}\label{tab2}
\end{table}

\begin{figure}[H]
	\centering
	\includegraphics[width=1\textwidth]{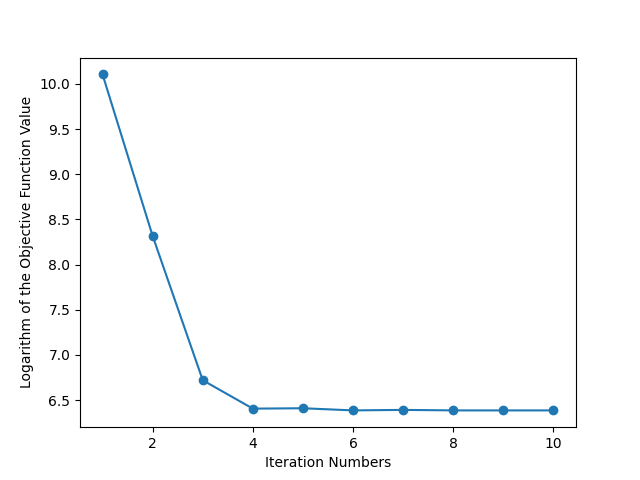}
	\caption{Variation of the logarithm of the objective function value $\ln Q(\boldsymbol{w})$ with the iteration numbers}
	\label{Fig.2}
\end{figure}

From the results in Table \ref{tab2}, it can be seen that the final expert weights obtained after iteration are $w_1=0.116$, $w_2=0.142$, $w_3=0.161$, $w_4=0.171$, $w_5=0.134$, $w_6=0.131$, and $w_7=0.145$, respectively. Meanwhile, the algorithm is iterated a total of 10 times in the process of obtaining the final solution that satisfies the error conditions. Among them, the change of the objective function value with the number of iterations is shown in Figure \ref{Fig.2}. It can be seen that the method converges faster, even in the face of a large amount of data, it can also meet the demand. When the number of experts, indicators and alternatives is large, the method can also get the results faster.

Based on the proposed optimization model, we simulated nearly 100 groups of data and conducted experiments on them respectively. Through the analysis of the experimental results we can find that: expert weights in descending order is consistent with the distances to the ``overall consistent score points'' in ascending order, that is, the larger the distance, the smaller the corresponding expert's weight, which is consistent with the idea when the model was established. The smaller the difference between the expert's scoring and the overall evaluation, the higher the similarity between the expert's evaluation and the consistent result, which means that the authority of its evaluation is higher, thus we give a high weight to this expert.

\section{Summary}
MAGDM is an important part of modern decision-making science, which is widely used in many fields. Determining expert weights objectively is a major difficulty, and the innovation of this paper is to provide a new method for determining expert weights by establishing an optimization model. We firstly determine the ``consistent score point'' which reflects the consistency of experts by using the evaluation of experts on each alternative in the model, and then use the sum of the distance between the expert scores and the points to construct the objective function which reflects the overall difference of the evaluation of experts, and the solution obtained by minimizing the objective function is the final expert weights. At the same time, we prove the existence and uniqueness of the solution of the model and give a brief introduction to the algorithm for solving the optimization model.

The expert weights obtained by the model have an important characteristic, i.e., the smaller the distance from the expert evaluation result to the ``overall consistent score point'', the larger the weights are, which is verified by numerical experiments in the paper. The weights have such a characteristic, avoiding the occurrence of individual experts' scoring maliciously deviating from the group evaluation in engineering practice. Meanwhile, in the real application, the user can change the constraints relatively freely according to the actual situation in order to enhance the adaptive and expressive ability of the model, which is another advantage of the optimization model.

The next research can be divided into two directions, one is that if the matrix $S$ does not satisfy the full rank of columns, we need to think about how to deal with the optimization model properly in practical engineering; the other is to explore how to apply the model to each specific alternative to improve the scientificity of decision-making and improve the existing MAGDM.

\begin{acknowledgements}
The authors would like to deliver thanks to Prof. Han Huilei and Prof. Huang Li for their assistance. This work is supported by the ``Key Research Program'' (Grant No.2022YFC3801300), the Ministry of Science and Technology, PRC.
\end{acknowledgements}

\bibliographystyle{spmpsci}

\end{document}